# SOME CASES OF THE FONTAINE-MAZUR CONJECTURE, II.

NIGEL BOSTON

Department of Mathematics, University of Illinois, Urbana, IL 61801

ABSTRACT. We prove more special cases of the Fontaine-Mazur conjecture regarding $p$-adic Galois representations unramified at $p$, and we present evidence for and consequences of a generalization of it.

## 0. Introduction.

Answering a longstanding question of Furtwängler, mentioned as early as 1926 [10], Golod and Shafarevich showed in 1964 [8] that there exists a number field with an infinite, everywhere unramified pro-$p$ extension. In fact it is easy to obtain many examples [8], [13], [23]. Very little is known, however, regarding the structure of the Galois group of such an extension.

In [7] Fontaine and Mazur conjecture (as a special case of a vast principle) that this Galois group can never be an infinite, analytic pro-$p$ group (i.e. linear over $\mathbf{Z}_p$). The idea is that a counter-example would produce an everywhere unramified Galois representation with infinite image, something that could not "come from algebraic geometry" (specifically the Galois action on a subquotient of an étale cohomology group, possibly with a Tate twist). Evidence for this conjecture has been published in [1] and [9].

As noted in [1], since infinite, analytic pro-$p$ groups contain a subgroup of finite index which is uniform [5, p.194] (in [1] the older terminology "$p$-saturated with integer values" was used), their conjecture is equivalent to the following.

**Conjecture 1.** (Fontaine, Mazur) There do not exist a number field $K$ and an infinite everywhere unramified Galois pro-$p$ extension $L$ such that $\text{Gal}(L/K)$ is uniform.

*Remarks.* The reduction of the problem from analytic to uniform pro-$p$ groups is useful in that uniform pro-$p$ groups have a simple, internal characterization, which analytic pro-$p$ groups lack. This is what is exploited in the main theorem below.

---

Partially supported by the Sloan Foundation and NSF grant DMS 96-22590. I thank God for leading me to these results.

Typeset by $\mathcal{A}_{\mathcal{M}}\mathcal{S}$-TEX

In [1], it was shown that Conjecture 1 is true if $K$ is of prime degree $\neq p$ over a subfield $F$ whose class number is prime to $p$ and if $L$ is Galois over $F$. One point of this paper is to describe how stronger results on fixed-point-free automorphisms produce strengthenings of this main theorem of [1], in particular allowing a weakening of the condition that $[K:F]$ be prime. We show that the result is true if $K/F$ is cyclic of any degree $n$ not divisible by $p$ (keeping the other hypotheses on $F$ and $L$).

More generally, we introduce the class of self-similar groups containing the class of uniform groups and show that the main theorem of [1] carries over with "uniform" replaced by "self-similar" if $K/F$ is cyclic of degree $n$ and a condition $H(\mathrm{Gal}(L/K), n)$ holds. The condition $H(G, n)$ is fundamental in the theory of fixed-point-free automorphisms and is conjectured always to hold.

This then takes the method of [1] about as far as it will go. Since self-similar groups arise as linear groups over more general rings than $\mathbf{Z}_p$, the above result suggests a natural generalization of the conjecture of Fontaine and Mazur. We study consequences of this generalization and end the paper with some related results.

## 1. The Main Theorem.

Since the main theorem concerns self-similar groups and the condition $H(G, n)$, we must first introduce these concepts.

**Definition.** Let $G$ be a pro-$p$ group. We say that $G$ is *self-similar* if $G$ contains a filtration of open, characteristic subgroups $G = G_1 \geq G_2 \geq ...$ with all $G_i/G_{i+1}$ abelian such that $\cap G_i = \{1\}$, together with isomorphisms $\phi_i : G_{i-1}/G_i \to G_i/G_{i+1}$ that commute with every continuous automorphism of $G$. (Note that since the factors $G_i/G_{i+1}$ are all abelian, we need only consider outer automorphisms.)

*Examples.* (1) A pro-$p$ group $G$ is called uniformly powerful (or uniform) if the above definition holds with $\phi_i$ being the map $x \mapsto x^p$ [5, p.64]. An example of such a group is $\ker(SL_2(\mathbf{Z}_p) \to SL_2(\mathbf{F}_p))(p > 2)$. Indeed, every profinite group linear over $\mathbf{Z}_p$ contains an open subgroup that is uniform [5, p.194].

(2) Let $G$ be a $\Lambda$-perfect pro-$p$ group [15] ($\Lambda = \mathbf{F}_p[[T]]$), e.g. $G = \ker(SL_2(\Lambda) \to SL_2(\mathbf{F}_p))(p > 2)$. Then $G$ is self-similar with $\phi_i$ being a so-called $T$-map [2].

Note that one consequence of the definition is that self-similar pro-$p$ groups are infinite (since $|G/G_i| = |G/G_2|^{i-1} \to \infty$ as $i \to \infty$).

**Definition.** Let $G$ be a pro-$p$ group and $n$ a positive integer. We say that $H(G, n)$ holds if there is a function of $n$ that is an upper bound for the derived length of every finite quotient of $G$ that admits a fixed-point-free automorphism of order $n$.

*Remarks.* $H(G, n)$ is conjectured to hold for all $G$ and all $n$ [22]. It is known to

hold when

(i) $n$ is prime or $= 4$, (any $G$),

or (ii) (for any $n$) the rank of the finite quotient groups of $G$ is bounded [21] ((ii) holds if $G$ is uniform [5,p.54]).

The following theorem generalizes the main theorem of [1].

**Theorem 1.** Suppose $K$ is a number field containing a subfield $F$ such that $K/F$ is cyclic of degree $n$ prime to $p$ and such that $p$ does not divide the class number $h(F)$ of $F$. Then there is no everywhere unramified pro-$p$ extension $L$ of $K$, Galois over $F$, with self-similar Galois group $\text{Gal}(L/K) = G$ such that $H(G,n)$ holds.

*Remarks* (i) The advances from [1] include the fact that we do not require the extension of $K$ to be uniform nor the degree of $K/F$ to be prime.

(ii) Theorem 1 carries over immediately to the situation where $L/K$ is unramified outside a finite set of primes, not including those above $p$, with the one modification that $h(F)$ must be replaced by the appropriate generalized class number (order of ray class group).

*Proof of Theorem 1.* Suppose that a counter-example to the theorem exists. Note that by the profinite version of Schur-Zassenhaus, the extension

$$1 \to \text{Gal}(L/K) \to \text{Gal}(L/F) \to \text{Gal}(K/F) \to 1$$

splits. Let us denote by $\sigma$ an element of $\text{Gal}(L/F)$ that maps under this splitting to a generator of $\text{Gal}(K/F)$ and write $G = \text{Gal}(L/K)$ for short.

Since the subgroups $G_i$ are characteristic in $G$, $\sigma$ acts on $G/G_i$ (by conjugation). Suppose it acts, for all $i$, with no fixed point other than the identity. Then since $H(G,n)$ holds, there is a bound on the derived length of $G/G_i$ as $i \to \infty$. This is false, since the maximal unramified pro-$p$ extension of any fixed derived length is a finite extension by repeated use of finiteness of class numbers.

Thus there is an $i$ and a nontrivial element $x \in G/G_i$ such that $\sigma(x) = x$. By picking the minimal such $i$ we ensure that $x$ maps to the identity in $G/G_{i-1}$, i.e. that $x \in G_{i-1}/G_i$. Since the map $\phi_{i-1}$ is $\sigma$-equivariant, there is a nontrivial fixed point in $G_{i-2}/G_{i-1}$ too, contradicting the minimality of $i$ unless $i = 2$, when we already have a fixed point in the abelian quotient $G/G_2$ of $G$. This, however, yields an unramified $C_p$-extension of $F$, contradicting $(p, h(F)) = 1$.

*Examples.* (1) Let $K$ be a quadratic field and $p$ an odd prime. If $L/K$ is any unramified pro-$p$ extension with $L$ Galois over $\mathbf{Q}$, then $\text{Gal}(L/K)$ is not self-similar.

(2) Let $K = \mathbf{Q}(\zeta_p)$ be the $p$th cyclotomic field, $p$ an odd prime. If $L/K$ is any unramified pro-$p$ extension with $L$ Galois over $\mathbf{Q}$, then $\text{Gal}(L/K)$ is not self-similar. Note that, for instance, by [23] the Hilbert $p$-class tower of $K$ is known to be infinite for $p = 157$.

## 2. A Generalization of the Fontaine-Mazur Conjecture.

Let $K$ be a number field, $p$ a rational prime, and $S$ a finite set of primes of $K$ containing none above $p$. Let $G_{K,S}$ denote the Galois group over $K$ of a maximal extension unramified outside $S$. Partially inspired the previous section wherein results for uniform groups carry over to self-similar groups, we conjecture the following:

**Conjecture 2.** Every continuous homomorphism $G_{K,S} \to GL_n(R)$, where $R$ is a complete, Noetherian local ring with finite residue field of characteristic $p$, has finite image.

*Remarks.* (i) Such rings $R$ are always quotients of some $W(k)[[T_1, ..., T_m]]$, where $W(k)$ is the ring of infinite Witt vectors over a finite field $k$ of characteristic $p$.
(ii) The case $R = \mathbf{Z}_p$ is due to Fontaine and Mazur [7] and generalizes Conjecture 1.

The reasons for believing this conjecture extend from elegance to evidence. It will make the deformation theory of $G_{K,S}$ particularly simple, as follows.

**Corollary to Conjecture 2.** Let $\bar{\rho} : G_{K,S} \to GL_n(k)$ be a continuous homomorphism with $k$ a finite field of characteristic $p$. The deformation theory of $\bar{\rho}$ factors through a finite quotient of $G_{K,S}$ and so can be computed as in [4].

Another consequence of the conjecture is described below. First, we must define an important subclass of pro-$p$ groups and describe how the critical cases of the Fontaine-Mazur conjecture concern this subclass.

**Definition.** A *just-infinite* pro-$p$ group is an infinite pro-$p$ group with all its nontrivial closed normal subgroups being open (in other words, it has no proper, infinite quotient).

*Examples.* The simplest such group is $\mathbf{Z}_p$. Another example is $\ker(SL_2(\mathbf{Z}_p) \to SL_2(\mathbf{F}_p))$ [22].

*Remarks.* An application of Zorn's lemma shows that every infinite, finitely generated pro-$p$ group has a just-infinite quotient [12]. Since the class of pro-$p$ groups linear over $\mathbf{Z}_p$ is closed under the operation of taking quotients, to prove the conjecture of Fontaine and Mazur mentioned in (ii) above it suffices to prove the following conjecture (for all possible $K$ and $S$):

**Conjecture 1'.** Every just-infinite pro-$p$ quotient of $G_{K,S}$ is not linear over $\mathbf{Z}_p$.

*Remarks.* Note that in effect we are reduced to considering minimal counterexamples to the Fontaine-Mazur conjecture. This is very similar to the reduction of certain questions in finite group theory to the case of finite simple groups. The conjecture pin-points a subclass of critical cases that specify what needs to be done in order to establish the Fontaine-Mazur conjecture. Since all just-infinite pro-$p$ groups satisfy the inequality of Golod and Shafarevich [22], the above shows that use of the failure of the Golod-Shafarevich inequality as in [9] to establish cases of the Fontaine-Mazur conjecture actually does not address the critical issue here.

There is, in analogy to the classification of finite simple groups, a (presently incomplete) classification of just-infinite pro-$p$ groups.

**The Classification of Just-Infinite Pro-$p$ Groups** [12],[22].

Traditionally, just-infinite pro-$p$ groups have been placed into four classes.
I. Solvable ones (that are necessarily linear over $\mathbf{Z}_p$).
II. Nonsolvable ones that are linear over $\mathbf{Z}_p$.
III. Nonsolvable ones that are linear over $\mathbf{F}_p[[T]]$.
IV. The rest! This so far means groups of Nottingham-type, Fesenko-type, and Grigorchuk-type.

*Remarks.* The groups in class I are simply $p$-adic space groups. Those in classes II and III are partially known (at least up to commensurability) [12].

If $k$ is a finite field, then the Nottingham group $N_k$ consists of the automorphisms of $k[[T]]$ given by $T \mapsto T + a_2 T^2 + a_3 T^3 + ...(a_i \in k)$. A group of Nottingham-type is an open subgroup of some $N_k$. The Fesenko group $S_q$ is the subgroup of $N_{\mathbf{F}_p}$ consisting of automorphisms $T \mapsto T + a_1 T^{1+q} + a_2 T^{1+2q} + ...(a_i \in \mathbf{F}_p)$, where $q$ is a power of $p$. A group of Fesenko-type is an open subgroup of some $S_q$ [6].

Groups of Grigorchuk-type are particular subgroups of $W_p$ where $W_p$ is the pro-$p$ automorphism group of the $p$-ary tree. For example, if $G_1 = C_2$ and $G_n = G_{n-1} \wr C_2$, then $W_2 \cong \lim_{\leftarrow} G_n$. Grigorchuk- type groups (sometimes called "branch") have complicated, recursively defined presentations.

Grigorchuk has recently proved that *every* just-infinite pro-$p$ group either is branch or contains an open subgroup of the form $H \times ... \times H$ (finitely many factors), where $H$ is hereditarily just-infinite (i.e. every open subgroup of $H$ is just-infinite).

**Corollary of Conjecture 2.** Any just-infinite pro-$p$ quotient of $G_{K,S}$ must be of type IV.

*Proof.* The map to such a quotient of type I, II, or III is a Galois representation of the sort disallowed by Conjecture 2. In fact, we need only use Class Field Theory to show that type I quotients do not arise and the Fontaine-Mazur conjecture to show that type II quotients do not arise.

*Remarks.* In a future paper examples of explicit Galois groups of some huge finite 2-extensions (for instance, degree $\geq 2^{21}$) of $\mathbf{Q}$ ramified only at a few odd primes (e.g. $3, 5, 7$) will be given and compared with quotients of known groups of the various types above.

## 3. Complements.

The above theorem adds to the information we have regarding Galois groups of unramified pro-$p$ extensions. Consider, in particular, the special case where $K$ is an imaginary quadratic field and $p$ an odd prime. Let $G$ be the Galois group over $K$ of a maximal unramified pro-$p$ extension. Here is a summary of the known properties of $G$:

I. $G$ is topologically finitely presented with generator rank $d(G)$ equal to its relation rank $r(G)$ [20].

II. If $1 \to R \to F \to G \to 1$ is a minimal presentation of $G$ in the sense that $F$ is free pro-$p$ on $d(G)$ generators, then $R \subseteq F_3$ ($F_i$ giving the Zassenhaus filtration on $F$) [13]. This implies that if $d(G) \geq 3$, then $G$ is infinite [13], by [14] not analytic, and by [15] not $\Lambda$-analytic.

III. Every open (i.e. finite index) subgroup of $G$ has finite abelianization (by finiteness of class numbers).

IV. $G$ has an automorphism of order 2 (namely complex conjugation) with no nontrivial fixed point on its abelianization.

V. Let $Q = G/G''$. For every normal subgroup $H$ of $Q$ with cyclic $Q/H$, the order of Ker $V$ is $[Q : H]$, where $V : Q/Q' \to H/H'$ is the transfer map [17].

VI. $G$ has no self-similar quotient stable under complex conjugation (see first example after Theorem 1).

Conjecturally (i.e. if the Fontaine-Mazur conjecture is correct), $G$ also has the property that every open subgroup of $G$ has no infinite analytic quotients. This generalizes III.

By Cohen-Lenstra heuristics [3], for any given odd prime $p$ and positive integer $d$, there should exist examples of such groups with $d(G) = d$. While this seems a natural class of pro-$p$ groups to attempt to classify, it is apparently not one familiar to pro-$p$ group theorists. Most classes considered satisfy the Golod-Shafarevich inequality, $r(G) > d(G)^2/4$ [8], or some refinement [13] of it, which our groups, as soon as $d(G) > 2$, do not. Indeed, one of the early ideas, due to Magnus [16], to show the nonexistence of infinite unramified pro-$p$ extensions, was to show the nonexistence of pro-$p$ groups satisfying III above. In [11], Itô produced pro-$p$ groups for which III holds, but they are of no use to us since they are uniform. (Note in fact that III holds for *all* nonabelian just-infinite pro-$p$ groups.) In a future paper David Perry will describe his work towards producing pro-$p$ groups satisfying I-VI and the conjectural generalization of III.

Finally, in [1], it was asked whether if $K$ is a number field with $p$-class field $L$ ($p$ odd) such that $p \mid h(L)$, there exists an unramified extension $M$ of degree $p$ of $L$ such that $M$ is Galois over $K$ and such that $\text{Gal}(M/K)$ has exponent equal to that of $\text{Gal}(L/K)$. It was noted there that this is true for $K$ quadratic and that the truth of the Fontaine-Mazur conjecture implies an affirmative answer to this question, when $K$ has an infinite $p$-class tower. Nomura [18] has shown this also holds if $p$ and $\ell$ are distinct odd primes such that the order of $p$ mod $\ell$ is odd and $K$ is an abelian $\ell$-extension of $\mathbf{Q}$.

As noted by Lemmermeyer, however, the answer to my question is in the nega-tive. He points out the example, due to Scholz and Taussky [19], of $\mathbf{Q}(\sqrt{-4027})$, whose 3-class tower terminates with Galois group isomorphic to the second of those listed in [1], namely $< x, y \mid y^{(x,y)} = y^{-2}, x^3 = y^3 >$. The point is that this group has a nonabelian subgroup of order 27 and exponent 9. Letting $K$ be the corresponding intermediate field, its 3-class field is an elementary abelian extension of degree 9 contained in no larger unramified extension with Galois group of exponent 3.